\documentclass[12pt]{article}
\usepackage{amssymb}
\usepackage{amsthm}
\usepackage{amsmath}
\usepackage{graphicx}
\usepackage{color}
\usepackage[T1]{fontenc}
\usepackage{amsfonts}
\usepackage{epsfig}
\usepackage{latexsym}
\usepackage{graphicx}
\usepackage{relsize}
\newtheorem{theorem}{Theorem}[section]
\newtheorem{lemma}[theorem]{Lemma}

\usepackage[parfill]{parskip}
\begin{document}
\title{\bf{A Revisit to the  Unit-Gompertz Distribution}}

\title{\textbf{Some More Properties of the Unit-Gompertz Distribution }}

\author{
        M. Z. Anis\thanks{Corresponding Author} \ \& Debsurya De\\
\footnotesize{\it{SQC \& OR Unit,  Indian Statistical Institute,}}\\
\footnotesize{\it{203, B. T. Road, Calcutta 700 108, India}}\\
\footnotesize{\it{E-mail:} {\rm zafar@isical.ac.in; debsurya01@outlook.com}}
}
\date{}
\maketitle

\begin{abstract}
In a recent paper, Mazucheli et al. (2019) introduced the unit-Gompertz (UG) distribution and studied some of its properties. It is a continuous distribution with bounded support, and hence may be useful for modelling life-time phenomena. We present counter-examples to point out some subtle errors in their work, and subsequently correct them. We also look at some other interesting properties of this new distribution. Further, we also study some important reliability measures and consider some stochastic orderings associated with this new distribution.
\end{abstract}

\section{Introduction}
During recent years there has been an increased interest in defining new generated classes of univariate continuous distributions. The works of Zografos and Balakrishnan (2009) and Risti\'{c}  and Balakrishnan (2012) may be mentioned as examples. Earlier, Eugene et al. (2002) introduced a general class of distributions generated from the logit of the beta random variable. The so called $ T-X $ transformation introduced by Alzaatreh et al. (2013) is another such attempt.\

In a similar vein, Mazucheli et al. (2019) introduced the unit-Gompertz (UG) distribution and studied some of its properties. More specifically, they considered the random variable $ X=e^{-Y}, $ where $ Y $ has the Gompertz distribution. They erroneously  claimed that its hazard rate function can admit all possible forms depending on the parameter. The support of this new distribution is $ (0, 1). $ It may be viewed as an alternative model for reliability studies where due to physical constraints such as design life of the system or limited power supply, distributions with a finite support might be required.  As an application, Jha et al. (2020) consider the problem of estimating multicomponent stress-strength reliability under progressive Type II censoring when stress and strength variables follow unit Gompertz distributions with common scale parameter. Jha et al. (2019) consider reliability estimation in a multicomponent stress –strength based on unit-Gompertz distribution. Kumar et al. (2019) are concerned with inference for the unit-Gompertz model based on record values and inter-record times.  \

However, some of the results presented in Mazucheli et al. (2019) are not entirely correct. We present counter-examples to point out these subtle errors in their work; and subsequently correct them in Section 2. We study conditional moments in Section 3. Other important   properties of this distribution are investigated in Section 4.  Reliability associated measures are studied in Section 5. Stochastic ordering are considered next in Section 6. Finally, Section 7 concludes the paper.

\

\section{Counter-examples and Corrections}	
For convenience, we shall stick to the notation of  Mazucheli et al. (2019). Let $ Y \left( \alpha, \beta\right)  $ be a non-negative random variable with Gompertz distribution having density function given by 
\[g\left( y\mid \alpha, \beta\right) = \alpha \beta\exp \left( \alpha+\beta y - \alpha e^{\beta y}\right),  \]
where $ y>0;  $ and $ \alpha > 0 $ and $ \beta >0 $ are shape and scale parameters, respectively. \

Using the transformation 
\[ X= e^{-Y}, \]
Mazucheli et al. (2019) obtained a new distribution with support on $ \left( 0, 1\right),  $ which they refer to as the \emph{unit-Gompertz distribution.} For completeness, we shall list down its pdf and cdf. The pdf of the unit-Gompertz distribution is given by 
\begin{equation}
f\left( x\mid\alpha, \beta\right) = \frac{\alpha\beta\exp \left[ -\alpha\left( 1/x^{\beta}-1\right) \right] }{x^{1+\beta}}; \mbox {~~~}\alpha>0, \beta>0, x \in (0, 1) \label{pdf}
\end{equation}
while its cdf is given by 
\begin{equation}
F\left( x\mid\alpha, \beta\right) =\exp \left[ -\alpha\left( 1/x^{\beta}-1\right) \right]; \label{cdf} 
\end{equation}	
and hence, the survival function is given by
\begin{equation}
\bar{F}\left( x\mid\alpha, \beta\right) =1-\exp \left[ -\alpha\left( 1/x^{\beta}-1\right) \right]; \label{sf} 
\end{equation}	
\subsection{Shape}

	 Mazucheli et al. (2019) erroneously stated (in their Proposition 1) that the pdf is log concave and unimodal over the \emph{entire}   support of $ X. $ As a counterexample, consider $ \alpha=0.25 ,$ and $ \beta=1. $ Let us consider the sign of the second derivative of $ \log f\left( x\mid \alpha, \beta\right)  $ at $ x=0.50. $ Routine calculation shows that $ \frac{d^{2}}{dx^{2}}\log f\left( x\mid \alpha=0.25, \beta=1\right)\mid_{x=0.50} >0, $ contradicting Proposition 1 and Equation (6) of Mazucheli et al. (2019). We correct their Proposition 1 as follows:
	\begin{theorem}
	The pdf of the UG distribution is log concave and unimodal for \\ $ x \in \left( 0,  \min {\left( \left( \alpha\beta\right)   ^{\frac{1}{\beta}}, 1\right) }\right]. $
		\end{theorem}
	\textbf{Proof:} The second derivative of $ \log f\left( x\mid \alpha, \beta\right)  $ is given by
	\[\frac{d^{2}}{dx^{2}}\log f\left( x\mid \alpha, \beta\right)=-\frac{\left( 1+\beta\right) }{x^{2}}\left( \frac{\alpha \beta}{x^{\beta}}-1\right) .\]
	Now observe that $ \alpha>0, $ $ \beta>0 $ and $ 0<x<1; $ hence$ \frac{\left( 1+\beta\right) }{x^{2}} $ is always $ >0. $\\ 
Hence, $\frac{d^{2}}{dx^{2}}\log f\left( x\mid \alpha, \beta\right)<0  $	if $ x\leq \left( \alpha\beta\right)   ^{\frac{1}{\beta}}. $
	
	This means that $ \log f\left( x\mid \alpha, \beta\right) $ is concave and unimodal for $ \alpha>0, $ $ \beta>0 $ and $ x \in \left( 0,  \min {\left( \left( \alpha\beta\right)   ^{\frac{1}{\beta}}, 1\right) }\right]. $ This completes the proof.\begin{flushright}
$\blacktriangleleft$
\end{flushright}
\begin{figure}[h]
\centering
\includegraphics[scale = 0.45]{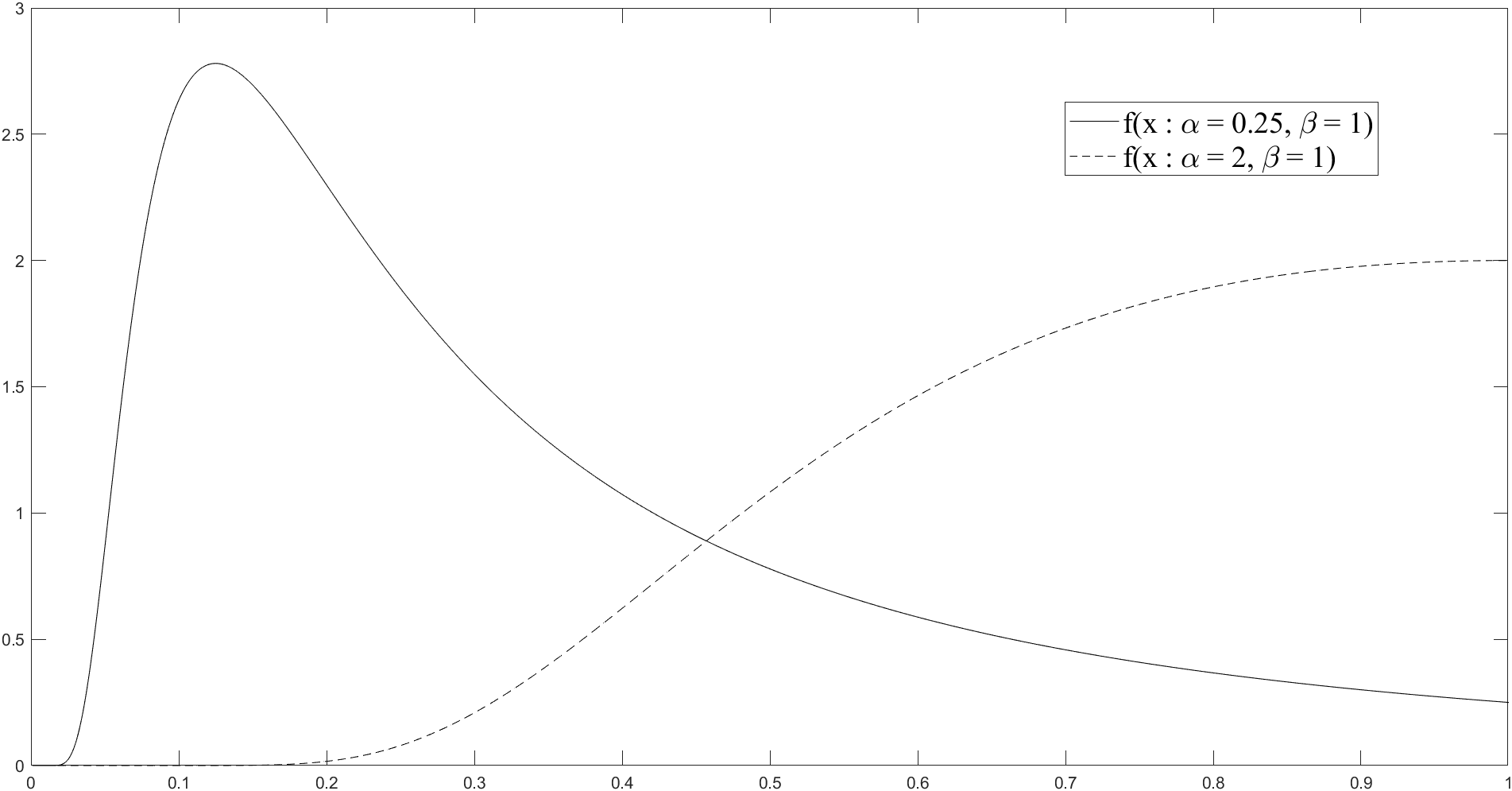}
\caption{Density Functions of unit-Gompertz distribution for $(\alpha = 0.25, \beta = 1)$ \& $(\alpha = 2, \beta = 1)$}
\end{figure}
Clearly, we see that the graph of  $ f(x;\alpha=2; \beta=1 ) $ is log-concave while the graph of $ f(x;\alpha=0.25; \beta=1 ) $ is not.

	\subsection{Mode}
	
	Consider the UG distribution with shape parameter $ \alpha = 3 $ and scale parameter $ \beta=1. $ According to equation (10) of Mazucheli et al. (2019), the modal point is $$ x_{0}= \left( \frac{\alpha \beta}{1+\beta}\right) ^{\frac{1}{\beta}}, $$ which in this particular case simplifies to 1.5. However, the support of the UG distribution is $ \left( 0, 1\right).  $ The plot of the density function is given in Fig. 2 below.
	\begin{figure}[h]
\centering
\includegraphics[scale = 0.45]{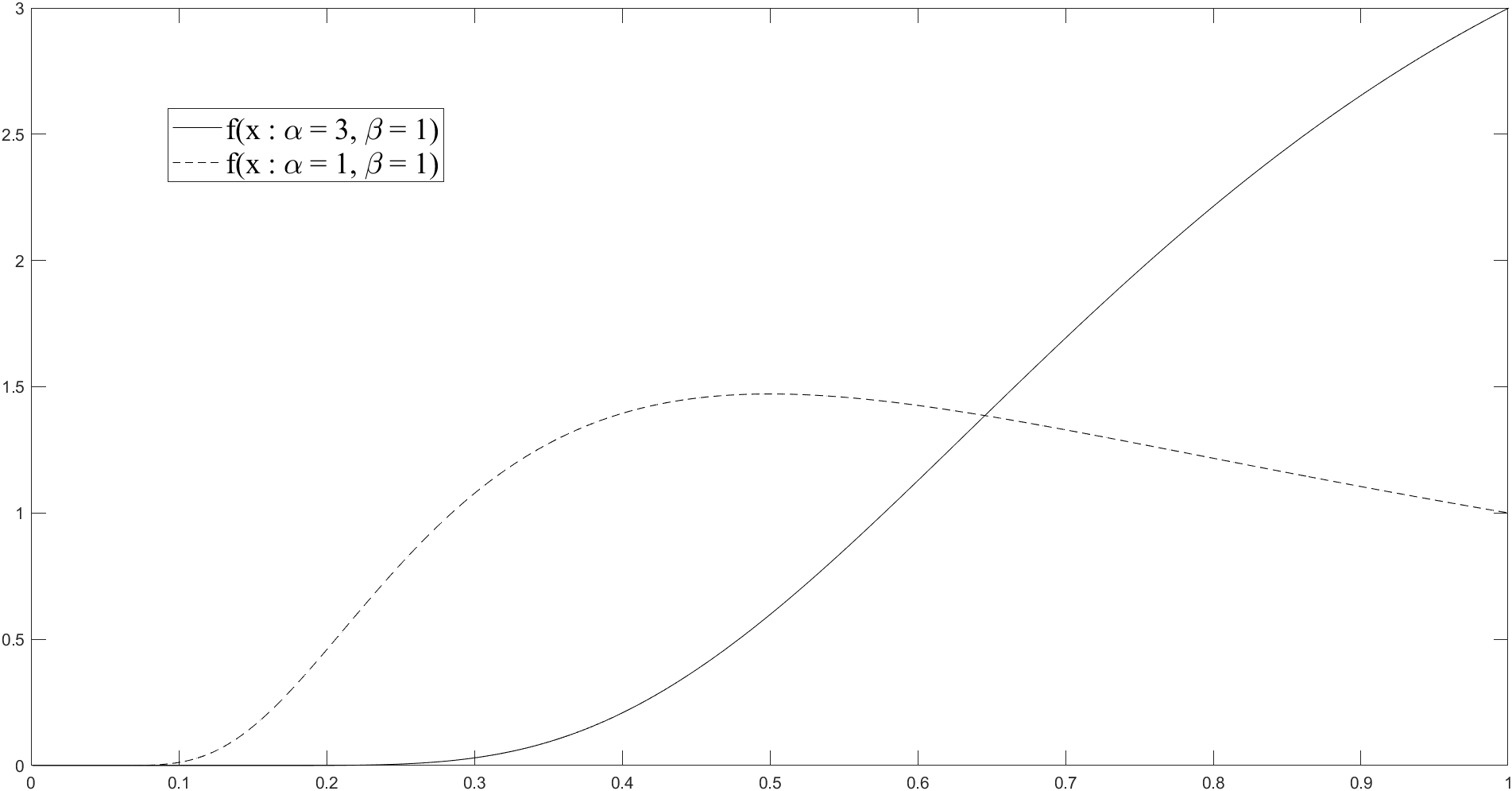}
\caption{Density Functions of unit-Gompertz distribution for $(\alpha = 3, \beta = 1)$ \& $(\alpha = 1, \beta = 1)$}
\end{figure}	\\
	It is easy to see from Fig. 2, that the mode of  $ f(x;\alpha=3; \beta=1 ) $  is at $ x=1,$ while the mode of $  f(x;\alpha=1; \beta=1 )$ is at $ x=0.5. $ We formalize this and now correct their result.\\
	 It is easy to see that the first derivative of  $ \log f\left( x\mid \alpha, \beta\right)  $ is given by
	\[\frac{d}{dx}\log f\left( x\mid \alpha, \beta\right)=-\frac{1+\beta}{x}+\frac{\alpha \beta}{x^{\beta+1}}.\]
	Hence, the mode of $ f\left( x\mid \alpha, \beta\right) $ is $ x^{\star}, $ the root of the equation 
	
	\begin{equation}
	\frac{d}{dx}\log f\left( x\mid \alpha, \beta\right) =0, \label{A10}
	\end{equation}
	
\noindent	
if $ x^{\star} \leq 1 $ where $ x^{\star}=\left( \frac{\alpha \beta}{1+\beta}\right) ^{\frac{1}{\beta}}. $ Hence the unique modal point is given by $$ x_{mode}= \min \left( x^{\star}, 1\right) . $$\

	\subsection{Hazard Rate  }
	
	Mazucheli et al. (2019) have erroneously mentioned (on page 27 of their paper) that  $ \lim_{x \rightarrow 1} h(x)= \alpha \beta;  $ and concluded that ``monotonically increasing shapes are possible for all values of $ \alpha>1 $ and $ \beta\geq1"; $ and ``possibly bathtub shapes of the hazard rate function will happen when $ \alpha\leq 0.5". $ Subsequently, they have used the result of Glaser (1980) to conclude (in their Theorem 3) that the hazard rate (HR) of the UG distribution is upside-down bathtub shaped.  They have sketched the hazard rate plot for different values of $ \alpha $  and $ \beta $ (in their Figure 2); but, it is important to note that not a single one of these graphs  seems to be upside-down bath-tub shaped. Fig. 3 shows the hazard plot of for a few selected values of  $ \alpha$ and $\beta.$\\

\begin{figure}[h]
\centering
\includegraphics[scale = 0.45]{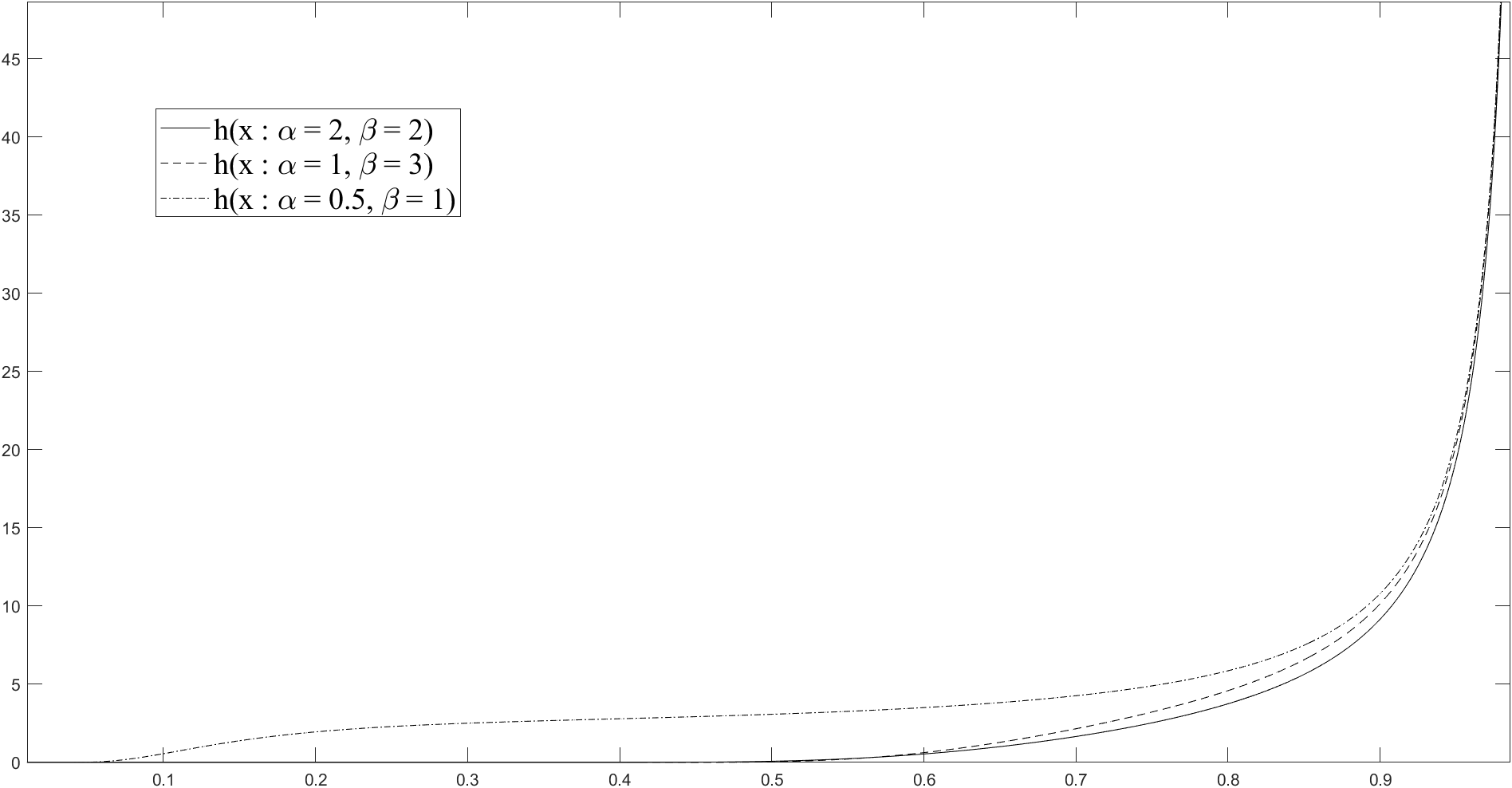}
\caption{Hazard Rate Functions of unit-Gompertz distribution for $(\alpha = 2, \beta = 2)$, $(\alpha = 1, \beta = 3)$ \& $(\alpha = 0.5, \beta = 1)$}
\end{figure}
We see from Fig. 3 that the hazard rate is \emph{not} upside-down bathtub shaped as enunciated in Theorem 3 of Mazucheli et al. (2019). This means that that their Theorem 3 is \emph{not}  correct.  In fact, the shape of the hazard rate for the UG distribution \emph{cannot} be obtained by appealing to Glaser's results. This is because Glaser's results are useful when the support of the distribution is $ \left( 0, \infty\right).$ However, in this case of the UG distribution,  the support is the finite interval $ \left( 0, 1\right).  $  Ghitany (2004) has obtained sufficient conditions to characterize the shape of the hazard rate when the support of the distribution is finite, say $ \left( 0, b\right).  $  But even Ghitany's theorem \emph{cannot} be applied to the UG distribution because $ f\left( 1\right) = \alpha \beta;  $ whereas Ghitany's theorem demands that $ f(b)=0, $ where $ b $ is the left-end support of $ f. $ \\ 
 
For the UG distribution, we have the hazard rate $$ h(x)= \frac{f(x)}{\bar{F}(x)}=\frac{\alpha\beta\exp \left[ -\alpha\left( 1/x^{\beta}-1\right) \right]}{x^{1+\beta}\left\lbrace 1-\exp \left[ -\alpha\left( 1/x^{\beta}-1\right) \right]\right\rbrace }. $$ Observe that $ \lim_{x \rightarrow 0+} h(x)=0 $ and $ \lim_{x \rightarrow 1-} h(x)= \infty.  $ Hence, there exists $ 0<M<1 $ such that $ h(x) $ is increasing in $ \left( M, 1\right),$ suggesting that the hazard rate cannot be upside-down bath-tub shaped.
	
\section{Conditional Moments}

The notions of conditional expectation (moment) and independence are routinely discussed in elementary probability and statistics courses at the undergraduate level. However, conditional moments are important in their own right, especially in probability theory and economics. Bryc (1996) considers conditional moment representations for dependent random variables. Domínguez and Lobato (2004) introduce simple and consistent estimation procedure for economic models directly based on the definition of conditional moments.  For the UG distribution with parameters $ \alpha $ and $ \beta, $ we have the following theorem:

\begin{theorem}
Let $ X $ follow the UG distribution with parameters $ \alpha $ and $ \beta. $ Then the conditional moment of $ X $ is given by
\[E\left( X^{n} \mid X>x\right) =\frac{e ^{\alpha}\alpha^{n/\beta}\left[ \Gamma\left( 1-\frac{n}{\beta}; \alpha \right) -\Gamma\left( 1-\frac{n}{\beta}; \frac{\alpha}{x^{\beta}} \right) \right], }{1-\exp \left[ -\alpha\left( 1/x^{\beta}-1\right) \right]}  \]
where $ \Gamma\left( s;x\right)  $ is the upper incomplete gamma function defined by 
\begin{equation}
	  \Gamma\left( s; x\right) =\int_{x}^{\infty}t^{s-1}e^{-t} dt.\label{inc-up-gamma}
	 \end{equation}

\end{theorem}

\textbf{Proof:} The conditional moment, $ E\left( X^{n} \mid X>t\right) $ can be written as
\[E\left( X^{n} \mid X>t\right) =\frac{1}{S\left( t\right) }I_{n}^{*}\left( t\right), \]
where
\[S\left( t\right)= 1- F\left( t\right)\]
and
\begin{eqnarray}
I_{n}^{*}\left( t\right) & = & \int_{t}^{1}y^{n} f(y) dy \label{eqa1} \\
&=& \alpha \beta e^{\alpha}\int_{t}^{1}y^{n-1-\beta}\exp \left( -\frac{\alpha}{y^{\beta}}\right)dy \nonumber \\
&=& \alpha  e^{\alpha}\int^{1/t^{\beta}}_{1}\frac{1}{z^{n/\beta}}e ^{-\alpha z} dz \nonumber \\
&=& e^{\alpha}\alpha^{n/\beta}\int^{\alpha/t^{\beta}}_{\alpha} u^{-n/\beta}e^{-u}du \nonumber \\
&=& e^{\alpha}\alpha^{n/\beta}\left[ \Gamma\left( 1-\frac{n}{\beta}; \alpha \right) - \Gamma\left( 1-\frac{n}{\beta};  \frac{\alpha}{t^{\beta}}\right)\right]; \label{eqn2} 
\end{eqnarray}
where $ \Gamma\left( s, x\right)  $ is the upper incomplete gamma function defined  in (\ref{inc-up-gamma}) above. This completes the proof.
\begin{flushright}
$\blacktriangleleft$
\end{flushright}
As applications of the concept of conditional moment, we may consider the evaluation of the mean residual life, the mean deviations about the mean and the median  and the expected inactivity time. These are discussed in the subsequent sections.

\section{Other Important Properties}

We shall now look at some other other distributional properties not considered in Mazucheli et. al (2019). Specifically we shall consider the mean deviation, the entropy,   the Lorenz,  Bonferroni and Zenga curves and order statistics.	
	
	\subsection{Mean Deviation}
	For empirical purposes, the shape of a distribution can be described by the so-called \emph{first incomplete moment,} defined by $ m_{1}\left( z\right) =\int_{0}^{z}xf\left( x\right) dx=\mu - I_{1}^{*}\left( z\right), $  where $ I_{n}^{*}\left( z\right) $ is defined in (\ref{eqa1}) above with $ n=1 $. This plays an important role in measuring inequality and is used to measure the dispersion and the spread in a
population from the center.  We shall first state a useful lemma.
	
\begin{lemma}
The mean deviation about any arbitrary point $ x_{0} $ is given by 
\[\delta\left( x_{0}\right)=E\left( \lvert X- x_{0}\rvert\right)= \int_{0}^{1}\mid x - x_{0}\mid f\left( x\right) dx = 2x_{0}F\left( x_{0}\right) -x_{0}F\left(0\right) -\mu +2 I_{1}^{*}\left( x_{0}\right) - x_{0}, \] 
where $ F\left(0\right) $ is defined as 0.
\end{lemma}	
The proof is simple and hence omitted.\\
	 Then the  mean deviation about the mean is given by

 $$\delta\left( \mu\right) = E\left( \lvert X- \mu \rvert\right)=2\mu F\left( \mu\right) -2 m_{1}\left( \mu\right) =2\mu F\left( \mu\right) -2\mu +2 I_{1}^{*}\left( \mu\right);  $$	 
	 and the  mean deviation 	about the median 
	$$ \delta\left( M\right) =E\left( \lvert X- M \rvert\right)=\mu F\left( \mu\right) -2 m_{1}\left( M\right)=2M F\left( M\right) -\mu +2 I_{1}^{*}\left( M\right) -M $$ 
	where $ \mu=E\left( X\right), $ $ M= Medain \left( X\right),   $ 
 $ m_{1}\left( z\right)= \int_{0}^{z}xf\left( x\right) dx  $ is the first incomplete moment;  and $ I_{1}^{*}(t) $ is as defined in  (\ref{eqa1}) above . The algebraic expressions for the mean and the median  have already been obtained by Mazucheli et al. (2019); and hence $ \delta\left( \mu\right)  $ and $ \delta\left( M\right)  $ can be easily evaluated numerically.

	\subsection{Entropies}

	An entropy is a measure of uncertainty of a random variable $ X. $ A large value of entropy implies greater uncertainty in the data. The concept of entropy is important in different subjects including communication theory, economics,   physics, probability and statistics. Nanda and Chowdhury (2020) provide a useful review. Several measures of entropy have been studied and compared in the literature. Two popular entropy measures are the Shannon and R{\'e}nyi entropies (Shannon (1951)and R{\'e}nyi (1961)). The R{\'e}nyi entropy of a random variable $ X $ with pdf $ f\left( \cdot\right)  $ is defined as 
	\[ I_{R}\left( \gamma\right) =\frac{1}{1-\gamma}\ln \int_{-\infty} ^{\infty}f^{\gamma}\left( x\right) dx,\]
	for $ \gamma>0 $ and $ \gamma \neq 1; $ while the Shannon entropy  is given by  $ E\left[ -\ln f\left( X\right) \right] .  $ It is a particular case of the R{\'e}nyi entropy for $ \gamma\uparrow 1. $\\

	 First, we shall calculate the R{\'e}nyi entropy. Towards this end, we compute
	 \[ \int_{0}^{1}\left[ f\left( x\right)\right] ^{\gamma} dx= \left( \alpha \beta e^{\alpha}\right) ^{\gamma}\frac{1}{\beta}\left( \alpha \gamma\right) ^{\frac{1}{\beta}\left[ 1-\gamma\left( 1+\beta\right) \right] }\Gamma\left(  \gamma+\frac{1}{\beta}\left(\gamma -1\right); \alpha\gamma\right)  \]
	 where $ \Gamma\left( s; x\right) $ 	 represents the upper incomplete gamma function defined in (\ref{inc-up-gamma}) above.\\
	 Then, the R{\'e}nyi entropy of $ X $ is given by 
	 \begin{eqnarray}
	 I_{R}\left( \gamma\right) & =& \frac{1}{1-\gamma}\left[\gamma \ln \left( \alpha \beta e^{\alpha}\right) +\ln \left( \Gamma\left( \gamma +\frac{1}{\beta}\left( \gamma-1\right) ;\alpha\gamma\right) \right)- \ln \beta +\frac{1}{\beta}\left\lbrace 1-\gamma\left( 1+\beta\right) \right\rbrace \ln (\alpha\gamma) \right]\nonumber \\
	 & = &\frac{1}{1-\gamma}\left[\alpha\gamma + \frac{\left( 1-\gamma\right)}{\beta}  \ln \alpha - \left( 1-\gamma\right) \ln \beta  +\frac{1}{\beta}\left\lbrace 1-\gamma\left( 1+\beta\right) \right\rbrace \ln (\gamma)\right. \nonumber\\
	 & & {}\mbox{~~~~~} +  \left. \ln \left( \Gamma\left( \gamma +\frac{1}{\beta}\left( \gamma-1\right) ;\alpha\gamma\right) \right) \right] \label{renyi}
	 \end{eqnarray}
	 
	Similarly, the Shannon entropy is given by 
	\begin{eqnarray*}
	E\left[ -\ln f\left( X\right) \right]&=& 1- \ln \left( \alpha \beta\right) - \left( 1+\beta\right) \frac{e^{\alpha}}{\beta}\Gamma\left( 0;\alpha\right).
	\end{eqnarray*}
	
	 The Shannon entropy can also be obtained by limiting $ \gamma \uparrow 1 $ in the R{\'e}nyi entropy  obtained above. 
	 
	 Song (2001) has shown that the gradient of the R{\'e}nyi entropy $ I_{R}^{\prime}\left( \gamma\right)=\left( d/d\gamma\right)I_{R}\left( \gamma\right)  $	is related to the log-likelihood by $ I_{R}^{\prime}\left( 1\right) =-\left( 1/2 \right) \mathrm{Var}\left[\left(  \log f\left( X\right)\right)  \right]. $ This equality and the fact that the quantity $-I_{R}^{\prime}\left( 1\right)  $ remains invariant under location and scale transformations motivated Song to propose $-2I_{R}^{\prime}\left( 1\right)  $ as a measure of the shape of a distribution. Taking the first derivative of (\ref{renyi}) and then limiting $ \gamma\uparrow 1 $ using L'Hospital's rule, one gets the expression
\begin{eqnarray*}
I_{R}^{\prime}\left( 1\right) &= & \frac{\beta +2}{2\beta}-\frac{1}{2}\left[ \alpha\left\lbrace \alpha -2\left( 1+\frac{1}{\beta}\right) \ln \alpha\right\rbrace - \left\lbrace \left[ \ln \alpha +e^{\alpha}\Gamma\left( 0; \alpha\right) \right] \left( 1+\frac{1}{\beta}\right) -\alpha \right\rbrace ^{2} \right.\\
& &  \mbox{~~~~~~~~~~~~~~~~~~~~}+ \left. e^{\alpha}\left( 1+\frac{1}{\beta}\right)^{2}\int_{\alpha}^{\infty}e^{-t}\left( \ln t\right) ^{2}dt \right] 
\end{eqnarray*}	
for the measure proposed by Song (2001). This measure plays a similar role as the kurtosis measure in comparing the shapes of various densities and measuring heaviness of tails.

\subsection{  Lorenz,   Bonferroni and Zenga curves}	

The Lorenz curve, introduced by Lorenz (1905), was proposed to measure the concentration of wealth. However, since then it has been used in many other areas. See, for example, Aaberge (2000), Jacobson et al. (2005) and  Groves-Kirkby et al. (2009) for its diverse use. In the field of reliability mention may be made of the works of Chandra and Singpurwalla (1984), Klefsj\"{o} (1984) and   Pham and Turkkan (1994).

Similarly, Bonferroni (1930) proposed a curve to measure wealth and income inequality. This curve also has applications in  demography, insurance and medicine. Giorgi and Crescenzi (2001) apply the   Bonferroni curve to analyze life-testing and reliability.

More recently Zenga (2007) introduced a new index of income inequality, $ Z(x), $ based on the ratio between the lower mean and the upper mean. $ Z(x) $ can be also interpreted as the difference in average age of components which has survived
beyond age $x$ from those which has failed before attaining age $x$, expressed in terms of average age of components exceeding age $ x. $ Hence, it can be viewed as a measure of proportional change in average age while switching over from survival before and after attaining age $ x. $ It is also related to the mean residual life function $ e_{F}\left( x\right)  $ as follows:

\[Z(x)=\frac{1}{F\left( x\right) } \left[ 1-\frac{E\left( X\right) }{x+e_{F} \left( x\right)}\right].\]
Nair and  Sreelakshmi (2012) discuss the Zenga curve in the context of reliability analysis.

These curves are defined as
\[L(p)= \frac{1}{\mu}\int_{0}^{q} xf(x) dx; \]

\[B(p)= \frac{1}{p\mu}\int_{0}^{q} xf(x) dx; \] and 
\[Z(x)= 1- \frac{\mu^{-}\left( x\right) }{\mu^{+}\left( x\right)}\]

respectively, where $ \mu^{+}\left( x\right)= E\left( X\mid X>x \right),  $  $ \mu^{-}\left( x\right)= E\left( X\mid X\leq x\right),   $    $q = F^{-1} ( p) $ and $ \mu=E(X). $ Then for the UG distribution, we have 
\begin{eqnarray*}
I(q) & =& \int_{0}^{q} x\frac{\alpha \beta}{x^{1+\beta}}\exp \left[-\alpha\left( \frac{1}{x^{\beta}}-1\right) \right] dx\\
& =& \alpha^{1/\beta}e^{\alpha}\int_{\alpha/q^{\beta}}^{\infty}e^{-z}z^{1-1/\beta-1}dz\\
&=& \alpha^{1/\beta}e^{\alpha}\Gamma\left( 1-\frac{1}{\beta}; \frac{\alpha}{q^{\beta}}\right),
\end{eqnarray*}
where $ \Gamma (s; x) $ is as defined in (\ref{inc-up-gamma}) above.
Hence, for $ \beta>1 $ we have \[L(p)=\frac{\alpha^{1/\beta}e^{\alpha}}{\mu}\Gamma\left( 1-\frac{1}{\beta}; \frac{\alpha}{q^{\beta}}\right);\]

\[B(p)=\frac{\alpha^{1/\beta}e^{\alpha}}{p\mu}\Gamma\left( 1-\frac{1}{\beta}; \frac{\alpha}{q^{\beta}}\right);\]
and 
\[Z(x)=1- \frac{\Gamma\left( 1-\frac{1}{\beta};\frac{\alpha}{x^{\beta}}\right) }{\left[\Gamma\left( 1-\frac{1}{\beta};\alpha\right) -  \Gamma\left( 1-\frac{1}{\beta};\frac{\alpha}{x^{\beta}}\right)\right] }\frac{\bar{F}(x)}{F(x)}. \]
\[Z(x)=1- \frac{\Gamma\left( 1-\frac{1}{\beta};\frac{\alpha}{x^{\beta}}\right) }{\left[\Gamma\left( 1-\frac{1}{\beta};\alpha\right) -  \Gamma\left( 1-\frac{1}{\beta};\frac{\alpha}{x^{\beta}}\right)\right] }\frac{\left\lbrace 1-\exp \left[ -\alpha\left( 1/x^{\beta}-1\right) \right]\right\rbrace }{\exp \left[ -\alpha\left( 1/x^{\beta}-1\right) \right]}. \]

\subsection{Order Statistics}

\newcommand{\bin}[2]
{
\left(
\begin{array}{@{}c@{}}
#1 \\ #2
\end{array}
\right)
}
It is well know that if $ X_{(1)} \leq X_{(2)}\leq \cdots \leq X_{(n)} $ denotes the order statistic of a random sample $ X_{1}, X_{2}, \cdots, X_{n} $ from a continuous population with cdf $ F_{X}(x) $ and pdf $ f_{X}(x), $ then the pdf of the $ j\mathrm{-th} $ order statistics is given by
\[f_{X_{(j)}}(x)=\frac{n!}{(j-1)!(n-j)!}f_{X}(x)\left[F_{X}(x)\right] ^{j-1} \left[1-F_{X}(x)\right] ^{n-j},\]
for $ j=1, 2, \cdots , n. $ Hence, the pdf of the $ j\mathrm{-th} $ order statistic  from the UG distribution will be given by 
\[f_{X_{(j)}}(x)=\frac{n!}{(j-1)!(n-j)!}\frac{\alpha \beta e^{\alpha j} e^{-\alpha j/x^{\beta}} } {x^{1+\beta}}\left[ 1- e^{\alpha}e^{-\alpha/x^{\beta}}\right] ^{n-j}.\]

The $ k\mathrm{-th} $ moment of $ X_{(j)} $ is obtained next. We have 
\begin{eqnarray*}
E\left( X^{k}_{\left( j\right) }\right)&=& \int_{0}^{1}x^{k}f_{X_{(j)}}(x) dx\\
&=& \int_{0}^{1}x^{k}\frac{n!}{(j-1)!(n-j)!}\frac{\alpha \beta e^{\alpha j} e^{-\alpha j/x^{\beta}} } {x^{1+\beta}}\left[ 1- e^{\alpha}e^{-\alpha/x^{\beta}}\right] ^{n-j} dx\\
&=&\frac{n!}{(j-1)!(n-j)!}\alpha e^{\alpha j}\int_{1}^{\infty}t^{-k/\beta}e^{-\alpha jt}\left[ 1-e^{\alpha}e^{-\alpha t}\right] ^{n-j}dt\\
&=&\frac{n!}{(j-1)!(n-j)!}\alpha e^{\alpha j}\int_{1}^{\infty}t^{-k/\beta}e^{-\alpha jt}\sum_{r=0}^{n-j} \bin {n-j}{r}\left( -e^{\alpha}e^{-\alpha t}\right) ^{r}dt\\
&=&\frac{n!}{(j-1)!(n-j)!}\alpha e^{\alpha j}\sum_{r=0}^{n-j} \bin {n-j}{r}\left( -1\right) ^{r}e^{r \alpha}\int_{1}^{\infty}t^{-k/\beta}e^{-\alpha t\left( j+r\right) }dt\\
&=&\frac{\alpha e^{\alpha j}n!}{(j-1)!(n-j)!}\sum_{r=0}^{n-j} \bin {n-j}{r}\left( -1\right) ^{r}e^{r \alpha}\left\lbrace \alpha\left( j+r\right) \right\rbrace ^{\frac{k}{\beta}-1}   \int_{\alpha\left( j+r\right) }^{\infty}y^{-k/\beta}e^{-y}dy\\
&=&\frac{\alpha e^{\alpha j}n!}{(j-1)!(n-j)!}\sum_{r=0}^{n-j} \bin {n-j}{r}\left( -1\right) ^{r}e^{r \alpha}\left\lbrace \alpha\left( j+r\right) \right\rbrace ^{\frac{k}{\beta}-1}  \Gamma\left( 1-\frac{k}{\beta}; \alpha\left( j+r\right) \right),  
\end{eqnarray*}
where $ \Gamma\left( s; x\right)  $	 is the upper incomplete gamma function defined in (\ref{inc-up-gamma}) above. Note that, the moments exists only when $ k< \beta. $
	\section{Some Other Important Reliability Functions}
	
	We shall now discuss the mean residual life (MRL), reversed hazard rate (RHR) and expected inactivity time (EIT) for the UG distribution.  The RHR and EIT may be looked as the dual properties of the HR and the MRL functions. We shall also investigate the monotonicity of these duals. Finally, we shall discuss stress strength reliability for the UG distribution.

\subsection{Mean Residual Life}	

An important ageing measure of interest is the mean residual life (MRL) function. It is defined as
\[e_{F}\left( t\right) = E\left( X-t\mid X >t\right). \]

Physically, it measures the expected remaining life time for a unit having already survived up to time $ t. $ It is easy to see that $ e_{F}\left( 0\right) =  E\left( X\right) =\mu,$ the mean of $ X. $ It can be calculated using the cdf or the pdf. Specifically, we have
\begin{eqnarray*}
e_{F}\left( t\right) & = & \frac{1}{\bar{F}\left( t\right)}\int_{t}^{\infty}\bar{F}\left( u\right)du \\ 
& = & \frac{\int_{t}^{\infty}uf\left( u\right) du}{\bar{F}\left( t\right) }-t \\ 
&=& \frac{1}{\bar{F}\left( t\right)}I_{1}^{*}\left( t\right) -t; 
\end{eqnarray*}
where $ I_{1}^{*}\left( t\right) $ can be obtained from (\ref{eqa1}) above with $ n=1.$ The expression of the MRL fucntion is quite complicated. However, since the hazard rate is increasing, it follows from Theorem 3 of Mi (1995) that the the MRL function $ e_{F}\left( t\right) $ is decreasing. Figure 4 shows the graph of the MRL function for some  combinations of $ \alpha$ and $\beta. $
\\
\begin{figure}[h]
\centering
\includegraphics[scale = 0.45]{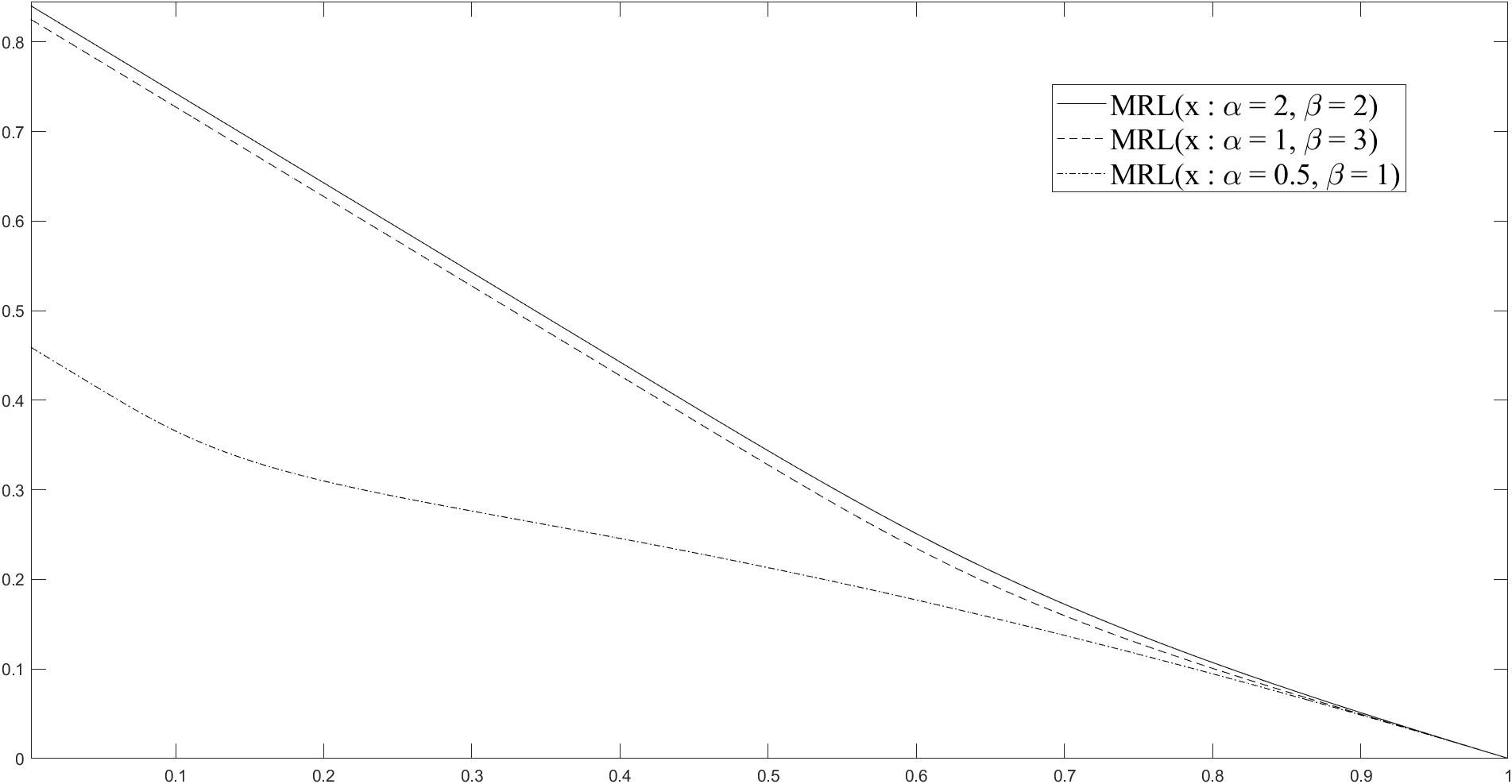}
\caption{Mean Residual Life Functions of unit-Gompertz distribution for \newline $(\alpha = 2, \beta = 2)$, $(\alpha = 1, \beta = 3)$ \& $(\alpha = 0.5, \beta = 1)$}
\end{figure}

	\subsection{Expected Inactivity Time}
	The Expected Inactivity Time (EIT) (also known as the mean past lifetime function) of a non-negative continuous random variable $ X $ with cumulative distribution function $ F\left( x\right)  $ is defined as 
	\[I\left( x\right) =E\left( x-X\mid X\leq x\right) = \frac{1}{F\left( x\right)}\int_{0}^{x}F\left( y\right)dy.\]
	Hence, $ I\left( x\right)  $ defines the mean waiting time for a device that failed in the interval $ \left[ 0, x\right] . $ In other words, this conditional random
variable shows the time elapsed from the failure of the component given that its
lifetime is less than or equal to $x.$ The EIT is a dual property of  the MRL; and is an important characteristic in many reliability applications. It has application in many disciplines such as survival analysis, actuarial studies and forensic science, to name but a few.  It is also of interest while describing different maintenance strategies. Chandra and Roy (2001) have shown that the EIT function cannot decrease on $ \left( 0, \infty\right); $ while Kundu and Nanda (2010) have studied some of its reliability properties. The non-parametric smooth estimation of the EIT function has been studied by Jayasinghe and Zeephongsekul (2013).\\

For the UG distribution, after detailed computation we

\[I\left( x\right) =\frac{e^{\alpha}\alpha^{1/\beta }}{\beta F\left( x\right)}\int_{\alpha /x^{\beta}}	^{\infty}e^{-u}\frac{du}{u^{1+1/\beta}}=\frac{e^{\alpha/x^{\beta}}\alpha^{1/\beta}}{\beta  }\Gamma\left( \frac{-1}{\beta};\frac{\alpha}{x^{\beta}}\right),\] 
which can be evaluated numerically.
	
	\subsection{Reversed Hazard Rate}
	The reversed hazard rate (RHR) of a non-negative continuous random variable $ X $ with pdf $ f\left( x\right)  $ and cdf $ F\left( x\right)  $ at time $ x $ is defined as 
	\[r\left( x\right) =\lim _{\Delta x\rightarrow 0}\frac{P\left( X>x-\Delta x\mid X\leq x\right) }{\Delta x}=\frac{f\left( x\right) }{F\left( x\right)}.\]
	Thus, $r\left( x\right)  $ defines the conditional probability of failure of a unit in $ \left( x-\Delta x,x\right)  $ given that the failure had occurred in in $ \left[ 0, x\right]. $ The RHR is a dual property of the HR. However, it should be noted that the trend in RHR is \emph{not} a direct indicator of the ageing pattern of a unit. The RHR has many interesting applications. Nanda and Shaked (2001) list the usefulness of the RHR while analyzing queuing systems. The RHR order arises naturally in economics and risk theory; see for example Eeckhoudt and Gollier (1995) and Veres-Ferrer and Pavía (2014). It is useful in estimating the survival function for left-censored lifetimes, see for example, Kalbfliesch and Lawless (1989). Irrespective of the shape of the hazard rate function, the RHR cannot increase on $ \left( 0, \infty\right) , $ as shown by Block et al. (1998). Testing the behaviour of the RHR is dealt with in Kayid et al. (2011).\\
	
	For the UG distribution, we have 
	\[r\left( x\right) =\frac{\alpha \beta}{x^{1+\beta}}.\]
	\subsection{Relationship}
	
	Chandra and Roy (2001) have proved the following:\\
	\begin{center}

		$ F\left( x\right)  $ is log-concave $ \Leftrightarrow $ $ r\left( x\right)  $ is decreasing $ \Rightarrow $ $I\left( x\right)  $ is increasing.
	\end{center}
	We shall use the above result to prove the following theorem:
\begin{theorem}
If $ X\backsim UG\left( \alpha, \beta\right) , $ then $ X $ has decreasing RHR (increasing EIT).
\end{theorem}	
	
\textbf{Proof:} We have
	\[\frac{d^{2}}{dx^{2}}\log F\left( x\mid \alpha, \beta\right)=-\frac{\left( 1+\beta\right) \alpha \beta}{x^{2+\beta}}<0\]
	since $ \alpha>0, $ $ \beta>0 $ and $ 0<x<1. $ Hence, $ F\left( x\mid \alpha, \beta\right) $ is log-concave. The theorem now follows from the result of Chandra and Roy (2001).
	\begin{flushright}
		$ \blacktriangleleft $
	\end{flushright}
	Figures 5 and 6 show respectively, the expected inactivity time and the reversed hazard rate for the UG distribution with parameters $ \alpha =0.25; 0.50; 0.75; 1.0  $ and $ \beta =1. $
\begin{figure}[h]
\centering
\includegraphics[scale = 0.45]{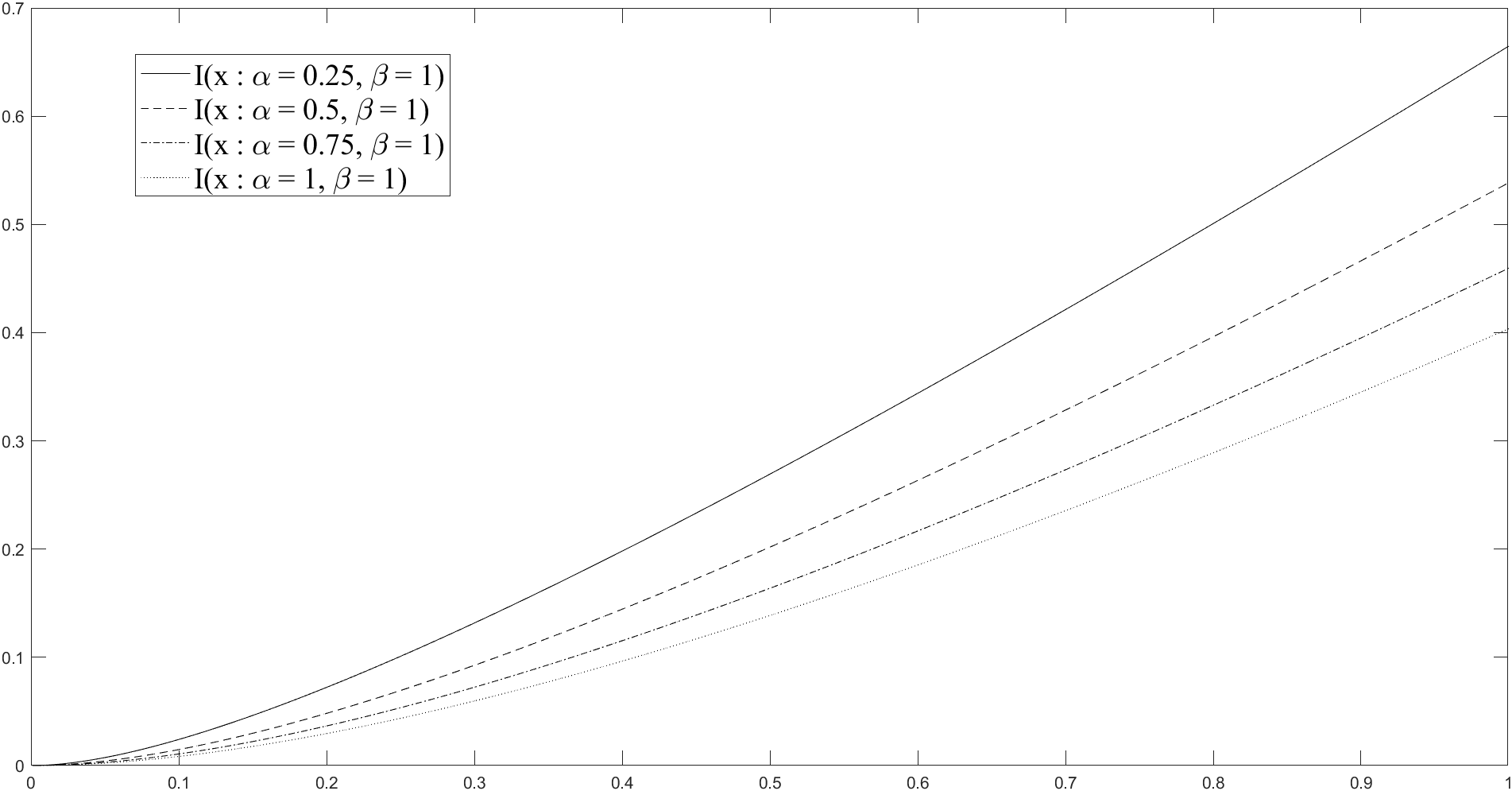}
\caption{Expected Inactivity Time Functions of unit-Gompertz distribution for \newline $(\alpha = 0.25, \beta = 1)$ , $(\alpha = 0.50, \beta = 1)$, $(\alpha = 0.75, \beta = 1)$ \& $(\alpha = 1, \beta = 1)$}
\end{figure}\\
\begin{figure}[h]
\centering
\includegraphics[scale = 0.45]{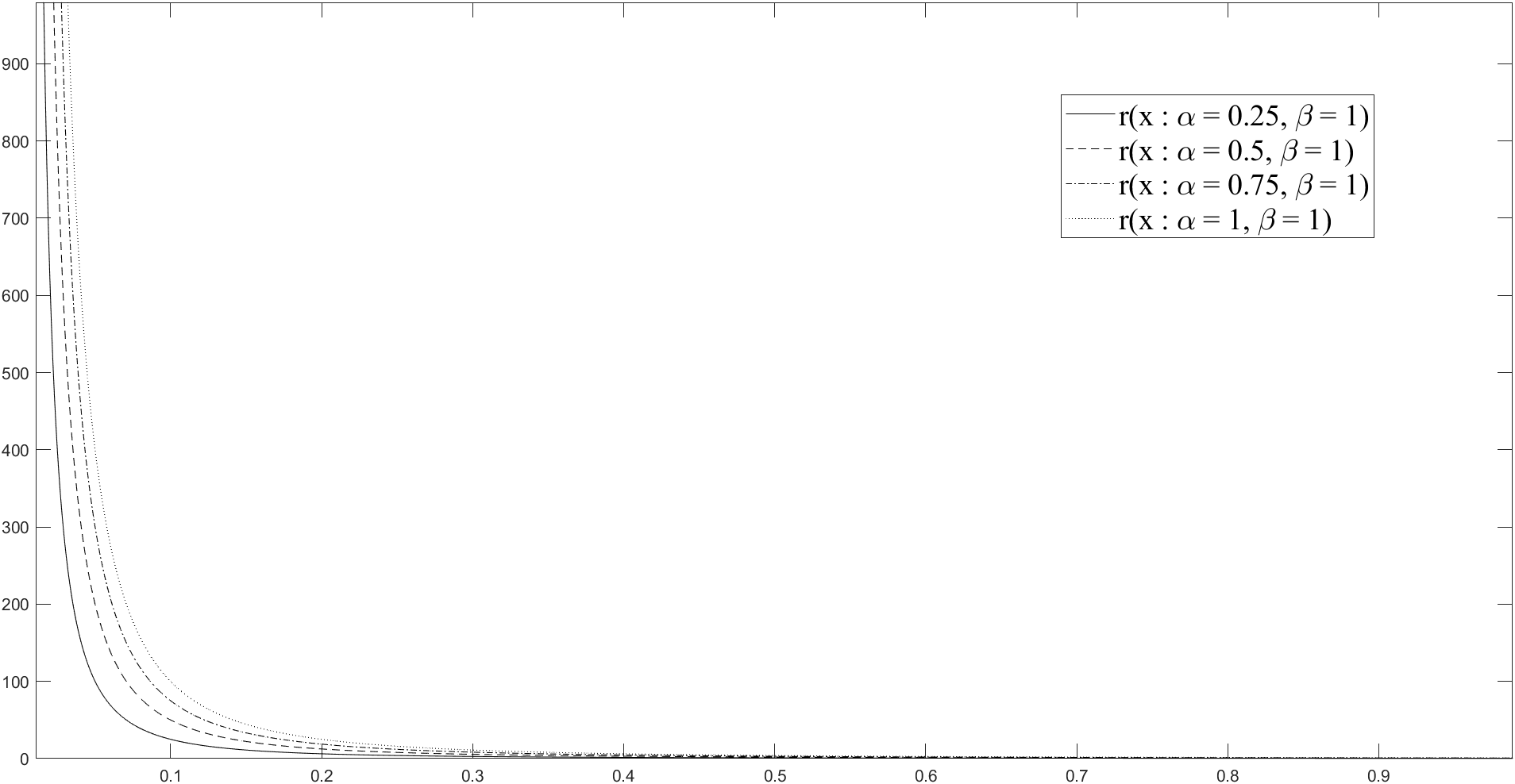}
\caption{Reverse Hazard Rate Functions of unit-Gompertz distribution for \newline $(\alpha = 0.25, \beta = 1)$ , $(\alpha = 0.50, \beta = 1)$, $(\alpha = 0.75, \beta = 1)$ \& $(\alpha = 1, \beta = 1)$}
\end{figure}
\newpage
\subsection{Stress Strength Reliability}	

Next, we derive the reliability $ R=Prob \left( Y<X\right) , $ where $ X \backsim UG \left( \alpha_{1}, \beta_{1}\right)  $ and $ Y \backsim UG \left( \alpha_{2}, \beta_{2}\right)  $ are independent random variables with distribution functions $ F_{X} $ and $ G_{Y} $ respectively.  Notionally, we may think of $ X $ as the strength and $ Y $ as the stress. Probabilities of this form have many engineering applications. 
\begin{eqnarray*}
R & = & Prob \left( Y<X\right)\\
& =& \int_{0}^{1} G_{Y}\left( x\right)f_{X} \left( x\right)dx\\
& =&\int_{0} ^{1}\frac{\alpha_{1}\beta_{1}}{x^{1+\beta_{1}}}\exp \left[-\alpha_{1} \left( \frac{1}{x^{\beta_{1}}}-1\right) \right] \exp \left[-\alpha_{2} \left( \frac{1}{x^{\beta_{2}}}-1\right) \right]dx\\ 
& =&\alpha_{1}\beta_{1}e^{\alpha_{1}+\alpha_{2}}\int_{0} ^{1}\frac{1}{x^{1+\beta_{1}}}\exp \left[- \left( \frac{\alpha_{1}}{x^{\beta_{1}}}+\frac{\alpha_{2}}{x^{\beta_{2}}}\right) \right]dx,
\end{eqnarray*}
which can be evaluated numerically. 
However, if the strength and stress distributions are independent random variables with with common scale parameter $ \beta, $ then we get a neat expression for $ R $ as $$ \frac{\alpha_{1}}{\alpha_{1}+\alpha_{2}} $$

	\section{Stochastic Orderings}
	Comparison of random variables based on their means, medians or variances is not very informative. The need to provide a more detailed comparison of two random quantities has been the origin of the theory of stochastic orders that has grown significantly during the last forty years. We shall begin by recalling some basic definitions. \\

	Let $X$ and $Y$ be random variables with distribution functions $F_{X}$ and $F_{Y},$ and survival
	functions $ \bar{F}_{X} $ and $ \bar{F}_{Y} $, respectively. Denote $ F^{-1}_{X}\left( u\right) = \sup \left\lbrace x: F_{X}\left( x\right)  \leq u \right\rbrace   $ and  $ F^{-1}_{Y}\left( u\right) = \sup \left\lbrace x: F_{Y}\left( x\right)  \leq u \right\rbrace   $  for $ u \in \left[ 0, 1\right] $ the right continuous inverses of $F_{X}$ and $F_{Y}.$\\

	\textbf{Definition: }A random variable $ X $ is said to be smaller than a random variable $ Y $ in the

	\begin{itemize}
		\item [(i)] usual stochastic order (denoted by $ X\leq_{st} Y $) if $ F_{X}\left( t\right)\geq  F_{Y}\left( t\right)  $ for all real $ t; $
		\item [(ii)] hazard rate order (denoted by $ X\leq_{hr} Y $) if $ \bar{F}_{X}\left( t\right)/ \bar{F}_{Y}\left( t\right) $ decreases in $ t; $ 
		\item [(iii)]  reversed hazard rate order (denoted by $ X\leq_{rh} Y $) if $ F_{X}\left( t\right)/ F_{Y}\left( t\right) $ decreases in $ t; $ 
		\item [(iv)]  mean residual life order (denoted by $ X\leq_{mrl} Y $) if $ \mu_{X}\left( t\right) \leq  \mu_{Y}\left(t\right); $
		\item [(v)]  expected inactivity time order (denoted by $ X\leq_{eit} Y $) if $ I_{X}\left( t\right) \geq  I_{Y}\left( t\right); $
		\item [(vi)] likelihood ratio order (denoted by $ X\leq_{lr} Y $) if $ f_{X}\left( t\right)/  f_{Y}\left( t\right) $ decreases in $t;$
		\item [(vii)] increasing convex order (denoted by $ X\leq_{icx} Y $) if $ \int_{x}^{+\infty}  \bar{F}_{X}\left( t\right)dt \leq \int_{x}^{+\infty}  \bar{F}_{Y}\left( t\right)dt, $ for all $ x,  $ provided the  two integrals exist;
		\item [(viii)] increasing concave order (denoted by $ X\leq_{icv} Y $) if $ \int_{-\infty}^{x}  F_{X}\left( t\right)dt \geq \int_{-\infty}^{x}  F_{Y}\left( t\right)dt, $ for all $ x,  $ provided the  two integrals exist;
		\item [(ix)] dispersive order (denoted by $ X\leq_{disp} Y $) if\mbox{ $ F_{X}^{-1} \left( \beta\right) - F_{X}^{-1} \left( \alpha \right) \leq  F_{Y}^{-1} \left( \beta\right) - F_{Y}^{-1} \left( \alpha \right) $} whenever $ 0 <\alpha \leq \beta <1; $
		\item [(x)]  stochastic-variability order (denoted as  $ X\leq_{st:icx} Y  $) if  $ X\leq_{st} Y $ and \mbox{$ \mathrm{Var}\left[ h\left( X\right) \right] \leq \mathrm{Var}\left[ h\left( Y\right) \right]  $} for any increasing convex function $h,$ provided the two variances exist;
		\item [(xi)] harmonic mean residual life order (denoted by  $ X\leq_{hmrl} Y $) if
		\[\left[ \frac{1}{x}\int_{0} ^{x}\frac{1}{m\left( u\right) }du\right] ^{-1}\leq  \left[ \frac{1}{x}\int_{0} ^{x}\frac{1}{l\left( u\right) }du\right] ^{-1} \mbox{~~~~ for all $x>0$}\]
		where $ m (u) $ and  $ l (u) $ are the mrl functions of the random variables $ X $ and $ Y $ respectively;
		\item [(xii)] star-shaped order (denoted by  $ X\leq_{ss} Y $) if $ E\left[ \phi\left( X\right) \right] \leq E\left[ \phi\left( Y\right) \right]  $ for all star-shaped functions $ \phi:\left[ 0, \infty\right) \longrightarrow \phi:\left[ 0, \infty\right), $ provided the expectations exist.
		\item [(xiii)] total time on test order (denoted by $ X\leq_{ttt} Y $) if and only if 
		\[\int_{0}^{F^{-1}_{X}(p)}\bar{F}_{X}(x) dx\leq \int_{0}^{F^{-1}_{Y}(p)}\bar{F}_{Y}(x) dx, \mbox{~~~~} p \in (0, 1).\]
	\end{itemize}
	The following implications are well known:
	\[X\leq_{eit}Y \stackrel {(e)}\Longleftarrow X\leq_{rh}Y\stackrel {(c)}\Longleftarrow X\leq_{lr}Y \stackrel {(a)}\Rightarrow  X\leq_{hr}Y \stackrel {(b)}\Rightarrow  X\leq_{mrl}Y \stackrel {(d)}\Rightarrow  X\leq_{hmrl}Y. \]
	\[ X\leq_{hr}Y \stackrel {(f)}\Rightarrow  X\leq_{st}Y \stackrel {(g)}\Rightarrow  X\leq_{ss}Y \stackrel {(h)}\Rightarrow  X\leq_{icx}Y\]
	\[X\leq_{st}Y  \stackrel {(i)}\Rightarrow  X\leq_{ttt}Y \stackrel {(j)}\Rightarrow  X\leq_{icv}Y\]
	
	The implications (a), (b) and (c) are given in Shaked and Shanthikumar (2007; page 43); implication (d) is in Shaked and Shanthikumar (2007; page 95) and implication (e) is given in Finkelstein (2002). The implications (f) is given in Shaked and Shanthikumar (2007; page 18) while the implications (g) and (h) can be found in Shaked and Shanthikumar (2007; page 205). The implications (i) and (j) can be found in Shaked and Shanthikumar (2007; page 224 and 225 respectively).\\
	
	The UG distributions are ordered with respect to the strongest likelihood ratio ordering as shown in the following theorem.\\
\begin{theorem}
Let $ X\backsim UG \left( \alpha_{1}, \beta\right) $ and $ Y\backsim UG \left( \alpha_{2}, \beta\right). $ If $ \alpha_{1} < \alpha_{2} ,$ then $$ X\leqq_{lr} Y \left( X\leqq_{hr} Y; X\leqq_{rh} Y; X\leqq_{mrl} Y; X\leqq_{eit} Y  \right) .  $$
	\end{theorem}	
	
	\textbf{Proof:} Let the corresponding pdfs be denoted by $ f_{X}\left( x\mid \alpha_{1}, \beta\right)  $ and $ f_{Y}\left( x\mid \alpha_{2}, \beta\right)  $ where $0<\alpha_{1} <  \alpha_{2} $ are the respective shape parameters and $ \beta>0 $ is the common scale parameter.\\
	Now observe that
	\[\frac{f_{X}\left( x\mid \alpha_{1}, \beta\right)  }{f_{Y}\left( x\mid \alpha_{2}, \beta\right)}=\frac{\alpha_{1}}{\alpha_{2}}\exp \left(\alpha_{1}-\alpha_{2}\right) \cdot \exp \left[ -x^{-\beta} \left(\alpha_{1}-\alpha_{2}\right)\right], \mbox{~~~} \alpha_{1}<\alpha_{2}. \]
	Hence, for $ \alpha_{1}<\alpha_{2}, $ we have
	\[\frac{d}{dx}\left(\frac{f_{X}\left( x\mid \alpha_{1}, \beta\right)  }{f_{Y}\left( x\mid \alpha_{2}, \beta\right)} \right)= \frac{\alpha_{1}}{\alpha_{2}}\exp \left(\alpha_{1}-\alpha_{2}\right) \cdot \beta \left(\alpha_{1}-\alpha_{2}\right) x^{-\beta - 1}  \exp \left[ -x^{-\beta} \left(\alpha_{1}-\alpha_{2}\right)\right]<0. \]
	This means $ {f_{X}\left( x\mid \alpha_{1}, \beta\right)  }/{f_{Y}\left( x\mid \alpha_{2}, \beta\right)} $ is decreasing in $ x. $ Hence, $ X\leq_{lr}Y .$\\
	
	The remaining statements follow from the implications given above. This proves the theorem.
	\begin{flushright}
		$\blacktriangleleft$
	\end{flushright}
	
	\section{Conclusion}
	
	This paper may be considered as an essentail follow-up paper of Mazucheli et al. (2019). It corrects some of the errors of the earlier paper. Some other important properties have been discussed. Since the proposed distribution can be used for modelling lifetime data, properties associated with lifetime distributions have been studied in the present work. It is hoped that this work will be a necessary complement to Mazucheli et al.(2019).

\end{document}